\documentclass[12pt]{amsart}

\usepackage[all,color]{xy}
\usepackage{pb-diagram}
\usepackage[mathscr]{eucal}

\usepackage{hyperref}
\hypersetup{
colorlinks=true, 
linktoc=all,     
linkcolor=blue,  
}
\usepackage{amsmath,amsfonts,amsthm,amssymb,amsthm}
\usepackage{xcolor}
\usepackage{pifont}
\usepackage{combelow} 

\usepackage[letterpaper,margin=1in,headheight=14pt,headsep=18pt]{geometry}
\usepackage[T1]{fontenc}
\usepackage{lmodern,xcolor}
\usepackage{amsmath,amssymb,amsthm,mathtools}
\usepackage{microtype}

\hypersetup{pdftitle={Concavity-Free L-infinity Estimates and Stability for Fully Nonlinear Complex Equations},pdfauthor={Ryan Mc Gowan}}

\theoremstyle{plain}
\newtheorem{example}{Example}[section]
\newtheorem{lemma}{Lemma}[section]
\newtheorem{corollary}{Corollary}[section]
\newtheorem{theorem}{Theorem}[section]
\newtheorem{proposition}{Proposition}[section]
\newtheorem{definition}{Definition}[section]

\numberwithin{equation}{section}

\DeclareMathOperator{\osc}{osc}

\DeclareMathOperator{\Vol}{Vol}

\newcommand{\cH}{i\partial\bar\partial}

\newcommand{\CP}{\mathbb {CP}}
\newcommand{\CC}{\mathcal C}
\newcommand{\calA}{\mathcal A}

\newcommand{\lam}{\lambda_\omega}

\allowdisplaybreaks[2]
\title[Concavity-free $L^\infty$ estimates and stability]{Concavity-Free $L^\infty$ Estimates and Stability for Fully Nonlinear Complex Equations}
\author{Ryan Mc Gowan}
\address{Department of Mathematics, Rutgers University, New Brunswick, NJ 08854, USA}
\email{rmm396@math.rutgers.edu}
\date{October 2026}

\begin{document}
\begin{abstract}
In this paper, we extend the $L^\infty$ estimates of Guo-Phong for fully non-linear equations on compact hermitian manifolds to a possibly indefinite reference form satisfying a uniform admissibility criterion. Furthermore, we also establish one-sided stability estimates. Neither result requires homogeneity or concavity of the nonlinear operator.  
\end{abstract}
\maketitle

\section{Introduction}

Since the development of De Giorgi-Nash-Moser theory for linear equations in divergence form, a priori estimates such as $L^\infty$-estimates have found themselves at the cornerstone of the theory of partial differential equations. Yau's seminal work \cite{Yau} resolving the Calabi conjecture highlighted their importance in complex geometry vis-\`a-vis the complex Monge-Amp\`ere equation on K\"ahler manifolds.  Ko{\l}odziej \cite{Kolodziej} then obtained a sharp version of Yau's $L^\infty$ estimates via pluripotential theoretic methods. Such estimates were extended to degenerate settings by Demailly-Pali \cite{DP} and Eyssidieux-Guedj-Zeriahi \cite{EGZ}, among others, and has become an important tool in the study of singular and degenerating K\"ahler metrics. These auxiliary Monge-Amp\`ere equations have found many applications in proving estimates in the hermitian setting; see, for instance, Guo-Song \cite{GS}. 

Recently, Guo-Phong-Tong \cite{GPTbounds,GPTstability} gave an alternate proof in a purely PDE theoretic manner, which applies to a large class of nonlinear equations on K\"ahler manifolds, opening the door to applying such techniques to problems in complex geometry. However, in the case of hermitian geometry, various analytical and geometric hurdles not present in the K\"ahler case add many technical difficulties. Previously, following Cherrier \cite{Cherrier}, Tosatti and Weinkove \cite{TW} established $L^\infty$-estimates in the hermitian setting. More recently, Guo and Phong \cite{GP} circumvented these difficulties by localizing the argument through the use of an auxiliary Monge-Amp\`ere equations to a ball with Dirichlet boundary conditions to then leverage the work of Caffarelli, Kohn, Nirenberg and Spruck \cite{CKNS}. 

The goal of this paper is to extend the estimates of \cite{GP} to possibly indefinite reference forms satisfying a uniform admissibility condition, while removing the homogeneity assumption on the nonlinear operator. General real reference  forms have previously been considered by Sz\'ekelyhidi \cite{Szekelyhidi} and Cheng and Xu \cite{CXdet,CXsub}, while Liu \cite{LiuABP} obtained related estimates. Our argument requires neither homogeneity nor concavity and also yields one-sided stability estimates. The main comparison inequality follows by integrating the operator along positive directions. We discuss the relation with the results of Cheng-Xu and Liu below.

We now state the main result of our paper. Let $(X,\omega)$ be a compact connected hermitian manifold without boundary and of complex dimension $n$.  Let $f:\Gamma\subset\mathbb R^n\to\mathbb{R}_+$ be a function satisfying the following three conditions:
\begin{enumerate}
\item[($1$)]
\makeatletter
\def\@currentlabel{$1$}
\makeatother
\label{cone} The set $\Gamma\subset \mathbb{R}^n$ is a symmetric convex cone with
\begin{equation*}
\Gamma_n\subseteq\Gamma\subseteq\Gamma_1,\end{equation*}
where $\Gamma_k$ is the cone of vectors $\lambda$ with $\sigma_j(\lambda)>0$ for $1\leq j\leq k$. Here $\sigma_j(\lambda)$ is the $j$-th elementary symmetric polynomial in $\lambda$.
\item[($2$)]
\makeatletter
\def\@currentlabel{$2$}
\makeatother
\label{symmetry} The function $f\in C^\infty (\Gamma)$ is symmetric and satisfies
\[
\frac{\partial f}{\partial \lambda_j}>0
\]
for each $j=1,\dots, n$ and $\lambda\in \Gamma$
\item[($3$)]
\makeatletter
\def\@currentlabel{$3$}
\makeatother
\label{positive} There is a constant $\gamma>0$ such that 
\[
\prod_{j=1}^n \frac{\partial f}{\partial \lambda_j}\geq \gamma
\]
for all $\lambda\in \Gamma.$
\end{enumerate}
In particular, we do not need assumption (2) in \cite{GP} for $f(\lambda)$, which we recall for the convenience of the reader and shall refer to henceforth as assumption ($2'$):
\begin{enumerate}
\item[($2'$)]
\makeatletter
\def\@currentlabel{$2'$}
\makeatother
\label{eq:structure} The function $f(\lambda)$ is symmetric in $\lambda\in \Gamma$ and is homogeneous of degree one. 
\end{enumerate}
As noted by Harvey and Lawson \cite{HL} the set of operators satisfying (1-3)+($2'$) is already quite large.

For a real $(1,1)$-form $A$, let $\lambda_\omega[A]$ denote the eigenvalues of $A$ with respect to the background metric $\omega$.

We fix a smooth real (possibly indefinite) $(1,1)$-form $\chi$ satisfying
\begin{equation}\label{eq:background}
\lambda_\omega[\chi-c_0\omega]\in\Gamma,
\qquad
\operatorname{Tr}_\omega(\chi)\leq K_\chi,
\end{equation}
for fixed $c_0>0$ and $K_\chi>0$. For $u\in C^2(X)$, consider
\[
\chi_u=\chi+i\partial\bar\partial u.
\]
We call $u$ admissible if $\lambda_\omega(\chi_u)\in\Gamma$ everywhere.  When an equation is imposed, we write its right-hand side as
\[
f(\lambda_\omega[\chi_u])=e^F.
\]
Throughout, $\omega^n=dV$, $V=\int_X dV$ and  $\|\cdot\|_p$ denotes the $L^p(X,\omega)$ norm for $1\leq p\leq \infty$.

Our first result extends the zero-order estimate of \cite{GP} to the indefinite background \eqref{eq:background}.

\begin{theorem}
\label{thm:Linfty}
Let $(X,\omega)$ be a compact hermitian manifold without boundary, and consider (\ref{eq:background}) where the operator $f(\lambda)$ satisfies conditions (1)-(3) above.  Let $p>1$ and $K>0$.  If $u\in C^2(X)$ is admissible, normalized to satisfy $\sup_X u=0$ and obeys 
\begin{equation*}
f(\lambda_\omega[\chi_u])=e^F,
\qquad
\|e^{nF}\|_p\leq K,
\end{equation*}
then one has the uniform estimate
\begin{equation*}
\|u\|_{L^\infty(X)}\leq C
\end{equation*}
for some constant $C>0$  depending only on $X$, $\omega$, $n$, $p$, $K$, $K_\chi$ and $\gamma$.
\end{theorem}

The second result is a stability estimate.  We first state it for smooth functions, echoing the result of Cheng and Xu \cite{CXdet}.  

\begin{theorem}
\label{thm:stability}
Assume the hypotheses of Theorem~\ref{thm:Linfty}.  Let $\varphi,\phi\in C^2(X)$ be admissible and satisfy $\sup_X\varphi=0$ and $\sup_X\phi\leq0$. If the following estimates hold:
\[
f(\lambda_\omega[\chi_\varphi])=e^G,
\qquad
\|e^{nG}\|_p\leq K,
\]
then for every
\begin{equation}\label{eq:alpha}
0<\alpha<\frac{p-1}{np+p-1},
\end{equation}
there is a constant $C$ depending only on the data above and on $\alpha$ such that
\begin{equation*}
\|(\phi-\varphi)_+\|_{L^\infty(X)}
\leq C\|(\phi-\varphi)_+\|_{L^1(X)}^\alpha.
\end{equation*}
\end{theorem}

A continuous weakly admissible version can also be deduced, which proves useful for compactness arguments.

\begin{definition}\label{def:weak}
A continuous function $v$ is called weakly $\chi$-admissible if, whenever $u\in C^2(X)$ touches $v$ from above at $x\in X$, one has $\lambda_\omega(\chi_u)\in\overline\Gamma$.
\end{definition}

\begin{corollary}\label{cor:weakstab}
In Theorem~\ref{thm:stability}, the comparison function $\phi$ may be replaced by any continuous weakly $\chi$-admissible function satisfying $\sup_X\phi\leq0$.
\end{corollary}

Finally, applying the theorem in both directions yields the usual two-sided estimate.

\begin{corollary}[Two-sided stability]\label{cor:twosided}
Let $\phi,\varphi\in C^2(X)$ be admissible and satisfy
\[
f(\lambda[\chi_\phi])=e^F,
\qquad
f(\lambda[\chi_\varphi])=e^G,
\qquad
\sup_X\phi=\sup_X\varphi=0,
\]
with
\[
\|e^{nF}\|_p+\|e^{nG}\|_p\leq 2K.
\]
Then, for every $\alpha$ satisfying \eqref{eq:alpha},
\[
\|\phi-\varphi\|_{L^\infty(X)}
\leq C\|\phi-\varphi\|_{L^1(X)}^\alpha.
\]
\end{corollary}

We discuss briefly the relation with recent work.  Cheng and Xu \cite{CXdet} proved a one-sided stability estimate of the same form, with the same range of exponents, under concavity, degree-one homogeneity and determinant majorization.  Thus the one-sided formulation and the exponent are not new. In a subsequent paper, Cheng and Xu \cite{CXsub} treated finite-slope Hessian quotient equations by using a strict viscosity subsolution, while retaining concavity. Liu \cite{LiuABP} obtained a complex Alexandrov-Bakelman-Pucci estimate and used it to derive related estimates for fully nonlinear equations with an indefinite admissible background.  In that application the operator is assumed concave, homogeneous of degree one and determinant-majorizing. The point we wish to emphasize is that the present paper circumvents the concavity condition by following \cite{GP}.

Guo and Song \cite{GS} introduced the sup-slope invariant and characterized solvability by a strict subsolution condition for a broad class of fully nonlinear elliptic equations.  Very recently, Guo and Song \cite{GSstability} used this strict sup-slope framework to obtain a stability estimate and then a direct gradient estimate on compact hermitian manifolds. The two settings are complementary: in \cite{GSstability} the strict sup-slope gap supplies the coercivity for the comparison, whereas here the global derivative-product bound \eqref{positive} provides that coercivity directly through the positive-ray inequality, with no equation-dependent subsolution.

It is useful to note that determinant majorization does not force concavity, even within the class covered by \eqref{symmetry}.  On $\Gamma_n$, set
\begin{equation*}
g(\lambda)=(\lambda_1\cdots\lambda_n)^{1/n},
\qquad
r(\lambda)=\Big(\sum_{j=1}^n\lambda_j^2\Big)^{1/2},
\qquad
f_\varepsilon=g+\varepsilon r.
\end{equation*}
Then $f_\varepsilon\geq g$, so determinant majorization holds.  Moreover
\[
(f_\varepsilon)_i
=\frac{g}{n\lambda_i}+\varepsilon\frac{\lambda_i}{r}
\geq\frac{g}{n\lambda_i},
\]
and hence
\begin{equation*}
\prod_{i=1}^n(f_\varepsilon)_i\geq n^{-n}.
\end{equation*}
At the point $(1,\ldots,1)$ and in every direction $v$ with $\sum v_i=0$,
\[
D^2f_\varepsilon[v,v]
=\left(-\frac1n+\frac{\varepsilon}{\sqrt n}\right)|v|^2.
\]
Thus $f_\varepsilon$ is not concave whenever $\varepsilon>1/\sqrt n$.  This gives a simple determinant-majorizing operator to which the results below apply but which is outside the concave framework of \cite{CXdet,LiuABP}.  We return to this example in Section~\ref{sec:examples}.

We now describe the proof.  If $A$ is admissible and $P>0$, integration along the positive ray gives
\begin{equation}\label{eq:rayintro}
f(\lambda_\omega[A+P])
\geq f(\lambda_\omega[A])
+n\gamma^{1/n}\det(\omega^{-1}P)^{1/n}.
\end{equation}
This is the only nonlinear comparison used in the paper.  Near a minimum of $u$, the condition \eqref{eq:background} allows us to subtract a small coordinate quadratic from $\chi$ and remain inside the cone. Comparison with a local Dirichlet Monge-Amp\`ere equation then gives Theorem~\ref{thm:Linfty}, exactly in the spirit of \cite{GP}.  For stability, one uses a global auxiliary equation whose right-hand side is concentrated on the set where $\phi-\varphi$ is large.  The perturbed form splits into an admissible base and a positive Monge-Amp\`ere increment; \eqref{eq:rayintro} then gives the contradiction.  The proof is carried out in Sections~2-4.\bigskip

The organization of this article is as follows. In section 2, we establish a useful linear algebraic inequality and recall the analytic estimates needed below. In section 3, we prove Theorem \ref{thm:Linfty} by localizing near a global minimum of the potential and comparing with an auxiliary Monge-Amp\`ere Dirichlet problem. In section 4, we then use Theorem \ref{thm:Linfty} to prove the stability estimate of Theorem \ref{thm:stability}. Finally, in section 5, we discuss some immediate consequences and examples highlighting the scope of these results when one drops the homogeneity assumption on the nonlinear operator.

\bigskip
\noindent\textbf{Acknowledgments.}
The author would like to thank their advisors, Professor Bin Guo and Professor Jian Song for introducing them to this problem, as well as their encouragement. This work forms part of the author's Ph.D. thesis at Rutgers University, and they are grateful to the Department of Mathematics for its generous support. \bigskip 

\noindent \textbf{AI Declaration:}  
The author is responsible for the original argument and takes full responsibility for the work presented here. AI was used to clean up the exposition and streamline the arguments.

\section{Preliminaries}\label{sec:prelim}

We begin with some linear algebraic arguments concerning the cone $\Gamma$. At a point $x\in X$, consider the spectral lift of $\Gamma$ defined by
\[
\CC_x=\{A=A^*: \lambda_{\omega(x)}(A)\in\Gamma\}
\]
A symmetric convex cone has a convex spectral lift. Further, $\CC_x$ is closed under the addition of positive semidefinite forms. In particular, given $A\in \overline{\CC_x}$ and $B\in \CC_x$, for each $0<t\leq 1$ one has
\begin{equation}\label{eq:mixedcone}
(1-t)A+tB\in\CC_x.
\end{equation}
Ellipticity in \eqref{symmetry} implies monotonicity under positive semidefinite increments.

We begin with a useful linear algebraic result to ensure \eqref{eq:rayintro} holds even in the abscence of concavity of the nonlinear operator $f$.

\begin{lemma}
\label{lem:ray}
For $x\in X$, and forms $A\in\CC_x$ and $P>0$, one has
\[
f(\lambda_\omega[A+P])
\geq f(\lambda_\omega[A])
+n\gamma^{1/n}\det(\omega^{-1}P)^{1/n}.
\]
\end{lemma}
\begin{proof}
For $t\in[0,1]$, consider the matrix $B(t)=A+tP$,
then $B(t)\in \mathcal{C}_x$ and $B'(t)=P>0$. For each fixed $t$, choose an $\omega$-unitary frame diagonalizing $B(t)$ 
\begin{align*}
\frac{d}{dt}f(\lambda_\omega[B(t)])&=\sum_{i=1}^n \frac{\partial f}{\partial \lambda_i}(\lambda_\omega[B(t)])\cdot p_{i\bar{i}}\\
&\geq n\left( \prod_{i=1}^n \frac{\partial f}{\partial \lambda_i}(\lambda_\omega[B(t)])\right)^{1/n}\left( \prod_{i=1}^n p_{i\bar{i}}\right)^{1/n}\\
&\geq n\gamma^{1/n}(\det(\omega^{-1}P))^{1/n}.
\end{align*}
The first inequality follows by the arithmetic-geometric mean inequality, and the second inequality follows by condition (\ref{positive}) and Hadamard's inequality. Integration in $t$ proves the result.
\end{proof}

The $L^1$-bound below on the admissible potential $u$ follows as a direct application of \cite[Lemma 8]{GP}.

\begin{lemma}\label{lem:L1}
For every admissible $u\in C^2(X)$ with $\sup_Xu=0$,
\begin{equation*}
\int_X(-u)\omega^n\leq C,
\end{equation*}
where $C$ depends only on $(X,\omega)$ and $K_{\chi}$.
\end{lemma}

We also recall some classical results which will prove helpful in the following sections:
\begin{itemize}
\item The classic theorem of Caffarelli, Kohn, Nirenberg and Spruck \cite{CKNS}, on a ball $B\Subset \mathbb{C}^n$ we have a unique smooth solution to the complex Monge-Amp\`ere Dirichlet problem
\begin{equation*}
\begin{cases}(\cH\psi)^n=g\omega^n,&\text{in }B\\
\psi=0,&\text{on }\partial B
\end{cases}
\end{equation*}
where $i\partial\overline\partial\psi\geq0$ and $g$ is both positive and smooth.  
\item By the work of Ko\l odziej \cite{Kolodziej}, we have the following exponential integral estimate: on a bounded pseudoconvex domain $\Omega\subseteq \mathbb{C}^n$ there are constants $\alpha$, $C$ depending only on $n$ and $\mathrm{diam}(\Omega)$ such that 
\begin{equation}\label{eq:exponentialinput}
\int_\Omega e^{-\alpha\psi}\omega^n\leq C
\end{equation}
for any $\psi\in \mathrm{PSH}(\Omega)\cap C^2(\Omega)\cap C^0(\overline\Omega)$ satisfying $\psi(z)=0$ for $z\in \partial \Omega$ and $\int_\Omega(\cH\psi)^n=1$.
\end{itemize}

Finally, a key component of the argument in the following sections is to study how badly mass may concentrate near global maxima/minima on $X$. The following abstract measure-theoretic proposition helps make this concentration more precise:

\begin{proposition}\label{lem:smallinc}
Let $m:[0,\infty)\to [0,\infty)$ be an nonincreasing function. Suppose there is some $T>0$ such that 
\[\begin{cases}
m(s)>0,&\text{for }0\leq s< T\\ 
m(s)=0,&\text{for }s\geq T.    
\end{cases}\] 
If for some $C,\delta>0$ and $0\leq s,t<T$    the following holds:
\[
t m(s+t)\leq Cm(s)^{1+\delta}
\]
then $T$ obeys the following bound:
\[
T\leq \frac{2Cm(0)^\delta}{1-2^{-\delta}}
\]
\end{proposition}

\begin{proof}
Starting with $s_0=0$, choose a sequence $\{s_j\}\subseteq [0,\infty)$ such that the subsequent terms obey
\[
s_{j+1}-s_j=2Cm(s_j)^{\delta}
\]
as long as $s_j<T$. Since $m(s_{j+1})\leq m(s_j)/2$ by a telescoping argument we find
\[
s_k = \sum_{j=0}^{k-1} s_{j+1}-s_{j} = 2C\sum_{j=0}^{k-1} m(s_{j})^\delta\leq 2Cm(0)^\delta \sum_{j=0}^{k-1} 2^{-\delta j}
\]
If $T\leq s_k$, at some step, the preceding geometric series argument still holds. If not, taking $k\to \infty$, we see
\begin{equation*}
s_\infty\leq
\frac{2Cm(0)^{\delta}}{1-2^{-\delta}}.
\end{equation*}
To finish, we show that $T\leq s_\infty$. Suppose by way of contradiction that $T>s_\infty$. As $m$ is nonincreasing, we have
\[
0<m(s_\infty)\leq m(s_j)\leq \frac{m(0)}{2^j}
\]
Hence, we see that $m(s_j)\to 0$ as $j\to \infty$, a contradiction. Hence $T\leq s_\infty$.  
\end{proof}

\section{Proof of Theorem \ref{thm:Linfty} via Auxiliary Equations}\label{sec:local}

Choose smooth positive functions
\[
\tau_j(t)=\frac{t+\sqrt{t^2+j^{-2}}}{2},
\]
so that $\tau_j\searrow t_+$ uniformly as $j\to \infty$. Let $v\in C^\infty(B)$, and define
\[
A_{j}(v)=\int_B\tau_j(v)e^{nF}\omega^n.
\]
Let $\psi_j$ be the solution to the auxiliary Dirichlet problem on a ball $B=B_1(0)$
\begin{equation*}
\begin{cases}
(\cH\psi_{j})^n
=\dfrac{\tau_j(v)}{A_{j}(v)}e^{nF}\omega^n,&\text{in }B 
\\ 
\psi_{j}=0,&\text{on }\partial B
\end{cases}
\end{equation*}
for some $\psi_j \in \mathrm{PSH}(B)$ provided to us by \cite{CKNS}. In particular,
\begin{equation*}
\int_B(\cH\psi_{j})^n=1,
\qquad \psi_{j}<0\quad\text{in }B, \qquad i\partial\overline\partial\psi_j\geq 0\quad \text{in }B.  
\end{equation*}

\begin{lemma}
\label{lem:comparison}
Let $\varepsilon$ be the constant defined by
\begin{equation*}
\varepsilon^{n+1}
=A_{j}(v)\,\gamma^{-1}\frac{(n+1)^n}{ n^{2n}}.
\end{equation*}
Let $v\in C^\infty(\overline B)$ and $v\leq 0$ on $\partial B$, and suppose that whenever $-i\partial \overline \partial v>0$ we have  
\begin{equation}\label{eq:detboundaux}
\det(\omega^{-1}(-\cH v))\leq\frac{e^{nF}}{n^n\gamma},
\end{equation}
then the following inequality holds throughout $\overline{B}$:
\begin{equation}\label{eq:comparison}
v_+\leq\varepsilon(-\psi_{j})^{\frac{n}{n+1}}.
\end{equation}
\end{lemma}
\begin{proof}
Consider the function $\Phi:B\to \mathbb{R}$ where
\[
\Phi=v-\varepsilon(-\psi_j)^{\frac{n}{n+1}}.
\]
Note that $\Phi\leq 0$ on $\partial B$. Suppose, by way of contradiction, it had a positive interior maximum within $B$. At that point, one finds after applying the chain rule
\begin{align*}
\cH v
&\leq\varepsilon\cH(-\psi_j)^{\frac{n}{n+1}}\\
&=-\frac{\varepsilon n}{n+1}(-\psi_j)^{-\frac{1}{n+1}}\cH\psi
-\frac{\varepsilon n}{(n+1)^2}(-\psi_j)^{-\frac{n+2}{n+1}}
i\partial\psi_j\wedge\bar\partial\psi_j\\
&\leq-\frac{\varepsilon n}{n+1}(-\psi_j)^{-\frac{1}{n+1}}\cH\psi_j.
\end{align*}
Consequently $-\cH v>0$ at that point, so \eqref{eq:detboundaux} applies.  Taking determinants, we see
\[
\left(\frac{\varepsilon n}{n+1}\right)^n(-\psi_j)^{-\frac{n}{n+1}}
\frac{\tau_j(v)}{A_{j}(v)}e^{nF}\leq \frac{e^{nF}}{n^n\gamma}.
\]
Since $e^{nF}>0$ and $\tau_j(v)\geq v$, we find
\[
v\leq A_{j}(v)\gamma^{-1}\frac{(n+1)^n}{n^{2n}\varepsilon^n}
(-\psi_j)^{
\frac{n}{n+1}
}
=\varepsilon(-\psi_j)^{
\frac{n}{n+1}
}
\]
contradicting the positivity of $\Phi$. Hence,  \eqref{eq:comparison} must hold.
\end{proof}

For $0\leq s< \mu$ where $\mu=\max_{\overline B}v$, consider the superlevel sets 
\[E(s)=\{x\in B:v(x)>s\}.\] 
We then define the weighted mass
\[
m(s)=\int_{E(s)}e^{nF}\omega^n.
\]
Note that $m(s)>0$ for $0\leq s<\mu$. We can extend $m(s)=0$ for $s\geq \mu$. We also define the weighted integral
\[
A(v,s)=\int_{E(s)}(v-s)e^{nF}\omega^n.
\]
We have that $A(v,s)$ is positive for $0\leq s<\mu$, and
\[A_j(v-s)\to A(v,s),\qquad \text{ as }j\to \infty.\]

\begin{proof}[Proof of Theorem \ref{thm:Linfty}.]

The idea is to first prove the result for $u\in C^\infty(X)$. The $C^2$ case then shall follow by approximation. 

To begin, let $M=-\min_Xu$ and consider a point $x_0\in X$ where this minimum is attained. Choose a holomorphic coordinate ball $B=B_1(x_0)$ and consider the Euclidean metric $\omega_{\mathrm{Euc.}}=i\partial \overline\partial|z|^2$. Consider some $0<\beta\ll1$ small enough to ensure
\begin{equation}\label{c0bound}
\beta \omega_{\mathrm{Euc.}}\leq c_0\omega    
\end{equation}
for some constant $c_0$. In particular, each of these choices can be chosen and scaled uniformly Then we note that
\begin{equation}\label{conetainment}
A=\chi-\beta\omega_{\mathrm{Euc.}}\in\CC_{x}.
\end{equation}
Define the potential
\begin{equation*}
v=-u-M+\beta(1-|z|^2).
\end{equation*}
We then see that $v$ obeys the properties
\begin{equation*}
v\leq0\ \text{on }\partial B,\qquad
\max_{\overline B}v=\beta>0.
\end{equation*}
Furthermore, whenever $P=-\cH v>0$ holds, we may apply Lemma~\ref{lem:ray} to the form $\chi_u=A+P$ to see that the following determinant bound holds:
\begin{equation*}
\det(\omega^{-1}(-\cH v))\leq\frac{e^{nF}}{n^n\gamma}.
\end{equation*}
Hence, we may apply Lemma \ref{lem:comparison} for each $0\leq s< \mu$ and the estimate \eqref{eq:exponentialinput} to get
\begin{equation*}
\int_{E(s)}\exp\left(
\alpha\frac{(v-s)^{\frac{n+1}{n}}}{A_j(v-s)^{1/n}}
\right)\omega^n\leq C,\qquad \alpha=C_1\frac{n^2}{n+1}\gamma^{1/n}.    
\end{equation*}
Applying Fatou's lemma, we see
\begin{equation}\label{eq:leveltrudinger}
\int_{E(s)}\exp\left(
\alpha\frac{(v-s)^{(n+1)/n}}{A(v,s)^{1/n}}
\right)\omega^n\leq C.
\end{equation}
Applying H\"older's inequality with exponents $p,p'$ gives
\begin{align}\label{areabound}
\frac{1}{A(v,s)^2}\int_{E(s)}(v-s)^{2n+2}e^{nF}\omega^n
\leq \left(\int_{E(s)} e^{pnF}\omega^n\right)^{1/p}\left(\int_{E(s)}\frac{(v-s)^{(2n+2)p'}}{A(v,s)^{2p'}}\omega^n\right)^{1/p'}.
\end{align}
On the right-hand side, the first integral is bounded since we assume $\|e^{nF}\|_p\leq K$ and the second integral is bounded by \eqref{eq:leveltrudinger}, hence we find
\begin{equation*}
\int_{E(s)}(v-s)^{2n+2}e^{nF}\omega^n\leq CA(v,s)^2.
\end{equation*}
Applying H\"older's inequality, we also have that \eqref{areabound} gives
\begin{align*}
A(v,s)
&\leq\left(\int_{E(s)}(v-s)^{2(n+1)}e^{nF}\omega^n\right)^{\frac{1}{2n+2}}
m(s)^{1-\frac{1}{2n+2}}\\
&\leq C A(v,s)^{\frac{1}{n+1}}m(s)^{1-\frac{1}{2n+2}}.
\end{align*}
Hence, we find
\begin{equation}\label{eq:powergain}
A(v,s)\leq Cm(s)^{1+\frac{1}{2n}}.
\end{equation}
If $t>0$, then we have the following relation between $m(s)$ and $A(v,s)$:
\begin{equation}\label{massarearel}
t\,m(s+t)\leq A(v,s).
\end{equation}
In particular, combining \eqref{eq:powergain} and \eqref{massarearel} becomes
\begin{equation*}
t\,m(s+t)\leq Cm(s)^{1+\frac{1}{2n}}.
\end{equation*}
We can thus apply Proposition \ref{lem:smallinc} with $T=\beta$ and $\delta=1/2n$ to get
\begin{equation}\label{eq:masslower}
m(0)\geq\left(\frac{\beta(1-2^{-\frac{1}{2n}})}{2C}\right)^{2n}=\nu>0.
\end{equation}
If $M>\beta$, then $-u>M-\beta$ on $E(0)$. 
By Lemma~\ref{lem:L1},
\[
\mathrm{Vol}_\omega(E(0)) \leq  \frac{1}{M-\beta}\int_{E(0)} (-u) \omega^n \leq\frac{C}{M-\beta}.
\]
H\"older's inequality and \eqref{eq:masslower} imply
\[
\nu
\leq m(0)\leq
K\left(\frac{C}{M-\beta}\right)^{1/p'}.
\]
Consequently,
\begin{equation*}
M\leq \beta+C\left(\frac{K}{\nu}\right)^{p'}.
\end{equation*}
The case $M\leq \beta$ is immediate. This proves Theorem~\ref{thm:Linfty} for $u\in C^\infty(X)$. 

For $C^2$ strictly admissible potentials, smooth approximation in the $C^2$ topology preserves admissibility for sufficiently large indices. The associated nonlinear densities converge uniformly, and the normalizations converge as well. Applying the smooth estimate with, for example, the common density bound $2K$ and passing to the limit gives the stated $C^2$ version. The same observation will apply to the controlled potential in the stability theorem.
\end{proof}

In particular, \eqref{eq:masslower} gives us the following corollary:

\begin{corollary}
\label{cor:noncollapse}
Fix $q>1$ and $K>0$. If $\psi$ is smooth, $\omega_\psi=\omega+\cH\psi>0$, and satisfies
\[
\omega_\psi^n = e^{nF}\omega^n, \qquad \|e^{nF}\|_q\leq K,
\]
then we have the uniform lower bound
\begin{equation*}
\int_X\omega_\psi^n\geq\nu>0
\end{equation*}
where the constant $\nu$ depends only on $(X,\omega)$, $q$, and $K$.
\end{corollary}
\begin{proof}
Let $\chi=\omega$ and nonlinear operator $f:\Gamma_n\to \mathbb{R}$ defined by
\[f(\lambda)=(\lambda_1\cdots\lambda_n)^{1/n}\]
For this operator $f$, we see 
\[\prod_{i=1}^n \frac{\partial f}{\partial \lambda_i}=n^{-n}>0.\]
We may choose $c_0=1/2$ in \eqref{c0bound} to ensure that $A\in\mathcal{C}_x$ in \eqref{conetainment}. To finish, apply \eqref{eq:masslower}.
\end{proof}

\section{Proof of Theorem \ref{thm:stability}}\label{sec:stabilityproof}

We first construct a normalized auxiliary Monge-Amp\`ere family with uniform bounds.

\begin{lemma}
\label{lem:globalaux}
Fix $q>1$. There are constants $B_q>0$ and $c_q>0$, depending only on $(X,\omega)$ and $q$, with the following property: If $\rho$ is smooth and satisfies
\begin{equation*}
\rho\geq\frac12,\qquad
\int_X\rho\omega^n=V,\qquad
\|\rho\|_q\leq2V^{1/q},
\end{equation*}
then the global hermitian Monge-Amp\`ere equation has a smooth solution $\psi$ satisfying
\begin{equation*}
(\omega+\cH\psi)^n=c\rho\,\omega^n,\qquad
\omega+\cH\psi>0
\quad\text{on }X,
\qquad \sup_X\psi=0.
\end{equation*}
such that 
\begin{equation*}
-B_q\leq\psi\leq0,\qquad c_q\leq c\leq2.
\end{equation*}
\end{lemma}
\begin{proof}
Smooth solvability, including the positive multiplicative constant $c$, follows from the work of Tosatti-Weinkove \cite{TW}.
At a maximum point of $\psi$, we have $\omega+\cH\psi\leq\omega$. Therefore $c\rho\leq1$ at that point, so $c\leq2$. It follows that
\[
\left\|\frac{(\omega+\cH\psi)^n}{\omega^n}\right\|_q
=\|c\rho\|_q\leq4V^{1/q}.
\]
Theorem~\ref{thm:Linfty}, applied to the Monge-Amp\`ere operator and the reference form $\omega$, provides the bound $B_q$. Further, Corollary~\ref{cor:noncollapse} gives
\[
cV=\int_X(\omega+\cH\psi)^n\geq\nu.
\]
We may therefore take
\[
c_q=\min\left\{1,\frac{\nu}{V}\right\}.\qedhere
\]
\end{proof}

\begin{proof}[Proof of Theorem \ref{thm:stability}]
As in the previous section, the goal will be to prove the statement for potentials $\varphi,\phi\in C^\infty(X)$ and then get $C^2$ regularity by an approximation argument. By Theorem~\ref{thm:Linfty}, we have the uniform estimate for some $U>0$ given by
\begin{equation*}
-U\leq\varphi\leq0.
\end{equation*}
Let $1<q<p$ and let $B_q,c_q$ be the constants from Lemma~\ref{lem:globalaux}. Choose $a>0$ depending on the background data to ensure $\chi-a\omega$ is admissible and fix
\begin{equation}\label{eq:fixedparameters}
L=1+U+aB_q.
\end{equation}
As before, consider the map
\begin{equation*}
v=\phi-\varphi,\qquad 
M=\sup_Xv.
\end{equation*}
and define 
\[
d=\int_Xv_+\omega^n.
\]
If $M\leq0$, the conclusion is immediate, so suppose that $M>0$. Since $\phi\leq0$, we see that $M\leq U$. Thus
\begin{equation}\label{eq:tchoice}
t=\frac{M}{8L}\quad\text{satisfies}\quad 0<t<\frac18.
\end{equation}
As before, consider the superlevel sets $E(M/4)$ and $E(M/2)$ of $v$. Choose $\eta\in C^\infty(X)$ with $0\leq\eta\leq1$, supported in $E(M/4)$, and equal to one on $\overline{E(M/2)}$. By Chebyshev's inequality,
\begin{equation*}
\Vol_\omega(E(M/4))\leq\frac{4d}{M}.
\end{equation*}
For notational convenience, set 
\begin{equation}\label{eq:Lambda}
\Lambda=\left(\frac{2}{a n t}\right)^n\frac{1}{\gamma c_q}.
\end{equation}
Suppose, by way of contradiction, that the following inequality holds:
\begin{equation}\label{eq:smallness}
\Lambda K\left(\frac{4d}{M}\right)^{\frac1q-\frac1p}
\leq\frac12V^{1/q}.
\end{equation}
Since $q<p$ we can apply H\"older's inequality to get
\begin{equation}\label{eq:densityq}
\|\Lambda\eta e^{nG}\|_q
\leq\Lambda K\Vol_\omega(E(M/4))^{\frac{1}{q}-\frac{1}{p}}
\leq\frac12V^{1/q}.
\end{equation}
Equation \eqref{eq:densityq} also gives
\[
\int_X\Lambda\eta e^{nG}\omega^n
\leq V^{1-1/q}\|\Lambda\eta e^{nG}\|_q\leq\frac V2.
\]
In order to apply Lemma \ref{lem:globalaux}, define
\begin{equation*}
\rho=1+\Lambda\eta e^{nG}-\frac{1}{V}\int_X \Lambda \eta e^{nG}\omega^n.
\end{equation*}
In particular, $\rho$ satisfies the conditions 
\begin{equation*}
\rho\geq\frac12,\qquad
\int_X\rho\omega^n=V,\qquad
\|\rho\|_q\leq\frac32V^{1/q}.
\end{equation*}
Hence, Lemma~\ref{lem:globalaux} gives the existence of $\psi,c$ and a positive form
\[P=\omega+\cH\psi\] 
satisfying the following properties 
\begin{equation}\label{eq:stabaux}
P^n=c\rho\,\omega^n,\qquad
-B_q\leq\psi\leq0,\qquad c_q\leq c\leq2.
\end{equation}
Consider the potential
\[
\phi_t=(1-t)\phi+ta\psi,
\]
Denote by $x_t$ the point which maximizes $\phi_t-\varphi$ and by $x_*$ the point which maximizes $v$. Consider the identity
\begin{equation}\label{eq:rewritemax}
\phi_t-\varphi=(1-t)v+t(a\psi-\varphi)
\end{equation}
Evaluating \eqref{eq:rewritemax} at $x_t$ and $x_*$ yields
\[
(1-t)v(x_t)\geq(1-t)M-t(L-1).
\]
Since $0<t<1/8$, it follows that
\begin{equation*}
v(x_t)\geq M-2tL=\frac{3M}{4}>\frac M2.
\end{equation*}
Thus $\eta(x_t)=1$ and $\rho(x_t)\geq\Lambda e^{nG(x_t)}$.
First suppose $\phi$ is $C^2$. At $x_t$,
\[
\chi_{\phi_t}\leq\chi_\varphi.
\]
We have the splitting
\begin{equation*}
\chi_{\phi_t}=A(t)+taP,
\qquad
A(t)=(1-t)\chi_\phi+t(\chi-a\omega).
\end{equation*}
Weak admissibility of $\phi$ and \eqref{eq:mixedcone} implies that $A(t)\in\CC_{x_t}$. In particular $\chi_{\phi_t}$ is admissible. Ellipticity from (\ref{positive}), Lemma~\ref{lem:ray}, and \eqref{eq:stabaux} combine to give
\begin{align}
e^{G(x_t)}
&\geq f(\lambda_{\omega(x_t)}[\chi_{\phi_t}])\nonumber\\
&\geq nta\gamma^{1/n}(\det(\omega^{-1}P))^{1/n}\nonumber\\
&\geq nta\gamma^{1/n}c_q^{1/n}\Lambda^{1/n}e^{G(x_t)}
=2e^{G(x_t)}.\label{eq:contradiction}
\end{align}
Here the final equality is the choice of $\Lambda$ in \eqref{eq:Lambda}. Since $e^G>0$, this is impossible.

The same argument applies to continuous weakly admissible $\phi$ without any regularization of $\phi$. Let
\[
k=\max_X(\phi_t-\varphi),\qquad
q_t=\frac{\varphi+k-ta\psi}{1-t}.
\]
Then $q_t$ is a $C^2$ upper test for $\phi$ at $x_t$. By Definition~\ref{def:weak},
\[
\chi+\cH q_t\in\overline{\CC_{x_t}}.
\]
At that point one has the exact identity
\begin{equation*}
\chi_\varphi
=(1-t)(\chi+\cH q_t)+t(\chi-a\omega)+taP.
\end{equation*}
The first two terms again form an interior admissible base. Lemma~\ref{lem:ray} now yields \eqref{eq:contradiction} directly. This proves that \eqref{eq:smallness} cannot hold for $M>0$, i.e. it must be the case that 
\[
\Lambda K\left(\frac{4d}{M}\right)^{\frac1q-\frac1p}>\frac{1}{2}V^{1/q}
\]
Substituting \eqref{eq:tchoice} and \eqref{eq:Lambda} into the previous equation gives
\[
\left(\frac{2}{an}\right)^n\frac{1}{\gamma c_q}K(8L)^n4^{\frac{1}{q}-\frac{1}{p}}\frac{d^{\frac{1}{q}-\frac{1}{p}}}{M^{n+\frac1q-\frac1p}}
>\frac12V^{1/q}.
\]
Consequently,
\begin{equation*}
M\leq C_qd^{\alpha(q)},\qquad
\alpha(q)=\alpha=\frac{p-q}{npq+p-q}.
\end{equation*}
This also holds when $d=0$, because continuity would otherwise give $d>0$. In particular, we see that as $q\downarrow1$, the exponent $\alpha(q)$ increases to
\[
\frac{p-1}{np+p-1}.
\]
This proves Theorem~\ref{thm:stability} for smooth controlled potentials, and the $C^2$ case follows like before.
\end{proof}

Finally, applying the theorem in both directions proves Corollary~\ref{cor:twosided}.

\section{Consequences and examples}\label{sec:examples}

\begin{corollary}\label{cor:compactness}
For fixed data as above, the family
\[
\calA_p(K)=\left\{u\in C^2(X):
\lam[\chi_u]\in\Gamma, \sup_Xu=0, \|f(\lambda_\omega[\chi_u])^n\|_p\leq K\right\}
\]
is relatively compact in $C^0(X)$. Its uniform limits are continuous and weakly $\chi$-admissible.
\end{corollary}
\begin{proof}
First, the family is relatively compact in $L^1$. To verify this on a hermitian manifold, let $L=\Delta_\omega$ and let $L^*$ be its formal adjoint with respect to $\omega^n=dV$. By Lemma~\ref{lem:L1}, 
\[\|u\|_{L^1(X)}\leq L_0,\] 
while a quick computation gives $(Lu)_-\leq K_\chi$. Hence
\begin{align*}
\int_X|Lu|dV
&=\int_XLudV+2\int_X(Lu)_-dV\\
&=\int_XuL^*1dV+2\int_X(Lu)_-dV
\leq L_0\|L^*1\|_\infty+2K_\chi V.
\end{align*}
The operator $1-L$ is invertible by the maximum principle and elliptic Fredholm theory. The derivative of its fixed Green kernel has the local bound $C\,\operatorname{dist}(x,y)^{1-2n}$. Integrating this bound in $x$ and applying Minkowski's inequality to $u=(1-L)^{-1}(u-Lu)$ gives
\[
\|u\|_{W^{1,s}(X)}\leq C(s)\|u-Lu\|_{L^1(X)}\leq C(s),
\qquad 1<s<\frac{2n}{2n-1}.
\]
Rellich compactness therefore gives an $L^1$-convergent subsequence.

Corollary~\ref{cor:twosided} makes this subsequence uniformly Cauchy, hence uniformly convergent to a continuous function.
\end{proof}

The compactness conclusion does not require convergence of the prescribed right-hand sides. 
Theorem~\ref{thm:stability} cannot in general be replaced by a two-sided estimate while keeping only the density bound for $\varphi$, as the following example shows:

\begin{example}\label{ex:onesided}
Let $X=\CP^1$, $\chi=\omega=\omega_{\mathrm{FS}}$, and take the one-dimensional operator $f(\lambda)=\lambda$. On the affine chart $U_0=\{z\in \mathbb{C}:z=Z_1/Z_0\}$, set, for $0<\varepsilon<e^{-1}$, the potential 
\[
v_\varepsilon(z)=-\frac1{2\log(\varepsilon)}
\log\left(\frac{|z|^2+\varepsilon^2}{1+|z|^2}\right).
\]
The argument of the logarithm is easily seen to be globally smooth and positive on $U_0$. Hence, $v_\varepsilon$ is smooth on $\CP^1$. Moreover,
\[
\omega_{\mathrm{FS}}+\cH v_\varepsilon
=\left(1+\frac1{2\log(\varepsilon)}\right)\omega_{\mathrm{FS}}
-\frac1{2\log(\varepsilon)}\cH\log(|z|^2+\varepsilon^2)>0.
\]
The second form is interpreted globally as the pullback of a Fubini-Study form by a projective automorphism. 
We have
\[
\sup_Xv_\varepsilon=0,\qquad
\inf_Xv_\varepsilon=-1,\qquad
\|v_\varepsilon\|_{L^1(X)}\longrightarrow0.
\]
The last assertion follows from $-1\leq v_\varepsilon\leq0$ and pointwise convergence to zero away from $z=0$. 
The fixed potential $\varphi=0$ has nonlinear density one. The positive part $(v_\varepsilon-\varphi)_+$ is zero, consistent with Theorem~\ref{thm:stability}, while
\[
\|v_\varepsilon-\varphi\|_{L^\infty(X)}=1
\quad\text{and}\quad
\|v_\varepsilon-\varphi\|_{L^1(X)}\longrightarrow0.
\]
\end{example}

The condition \ref{c0bound} is strictly weaker than positivity when $\Gamma$ is larger than $\Gamma_n$. For example, on a flat complex three-torus with $\omega=\omega_1+\omega_2+\omega_3$, take
\[
\chi=3\omega_1+3\omega_2-\omega_3.
\]
Its eigenvalue vector belongs to $\Gamma_2$, since
\[
\sigma_1(3,3,-1)=5>0,\qquad \sigma_2(3,3,-1)=3>0.
\]
Here, one can choose $c_0=1/4$: the shifted eigenvalues are $(11/4,11/4,-5/4)$, with positive first and second elementary symmetric functions. The form nevertheless has a negative eigenvalue everywhere.

For the normalized complex Hessian operator
\[
f(\lambda)=\left(\frac{\sigma_k(\lambda)}{\binom nk}\right)^{1/k}
\quad\text{on }\Gamma_k,
\]
the standard structural inequalities imply \eqref{eq:structure}; this class is included in the framework of \cite{GPTbounds,GP}. Writing $f(\lambda_\omega(\chi_u))=e^F$, the quantity $e^{nF}$ is a constant multiple of $\sigma_k(\lam(\chi_u))^{n/k}$. Thus the density hypothesis is a condition on the Hessian density in $L^{np/k}$, and it should not be confused with an $L^p$ bound on $\sigma_k$ with the same numerical exponent.

\begin{example}\label{ex:nonconcave}
For $n\geq 2$, on the positive cone $\Gamma_n$, let
\[
g(\lambda)=(\lambda_1\cdots\lambda_n)^{1/n},
\qquad
r(\lambda)=\left(\sum_{j=1}^n\lambda_j^2\right)^{1/2},
\qquad
f_\varepsilon(\lambda)=g(\lambda)+\varepsilon r(\lambda),
\]
where $\varepsilon>0$.  The function $f_\varepsilon$ is smooth, symmetric, positive and homogeneous of degree one, and
\[
f_\varepsilon(\lambda)\geq(\lambda_1\cdots\lambda_n)^{1/n}.
\]
Thus it satisfies determinant majorization on $\Gamma_n$. Furthermore,
\[
(f_\varepsilon)_i
=\frac{g}{n\lambda_i}+\varepsilon\frac{\lambda_i}{r}>0,
\]
so that
\[
\prod_{i=1}^n(f_\varepsilon)_i
\geq
\prod_{i=1}^n\frac{g}{n\lambda_i}
=\frac1{n^n}.
\]
Hence the structural condition \eqref{eq:structure} holds with $\gamma=n^{-n}$, independently of $\varepsilon$. The operator is nevertheless not concave for $\varepsilon>1/\sqrt n$.  Indeed, at $e=(1,\ldots,1)$ and for a nonzero vector $v$ satisfying $\sum_i v_i=0$,
\[
D^2g(e)[v,v]=-\frac1n|v|^2,
\qquad
D^2r(e)[v,v]=\frac1{\sqrt n}|v|^2.
\]
Therefore
\[
D^2f_\varepsilon(e)[v,v]
=\left(-\frac1n+\frac{\varepsilon}{\sqrt n}\right)|v|^2>0
\]
when $\varepsilon>1/\sqrt n$.  Thus determinant majorization, even together with the product condition, does not imply concavity.  This example is covered by Theorems~\ref{thm:Linfty} and \ref{thm:stability}, but not by results whose structural assumptions include concavity.
\end{example}

The proofs use an admissible base point before adding the positive auxiliary Hessian. A completely arbitrary indefinite $\chi$ outside $\Gamma$ does not provide such a base point. The present results therefore do not claim uniformity as the margin in \eqref{eq:background} collapses.

If $\chi$ itself is not admissible but there is a fixed smooth $q$ such that
$\widehat\chi=\chi+\cH q$ has the required margin, the equation can be rewritten with potential $u-q$ and reference form $\widehat\chi$. The uniform estimate then depends on that representative and on $\osc_Xq$. This observation does not supply a uniform admissible representative for a degenerating family.

\end{document}